\documentclass[12pt]{article}
\usepackage{amsmath, amssymb,amsfonts,bbm,verbatim,multirow}
\usepackage{epic,eepic,psfrag,epsfig}
\usepackage[sf,bf,SF,footnotesize]{subfigure}
\usepackage{graphicx}
\usepackage{amsthm}
\usepackage{amscd}
\usepackage{epsfig}
\usepackage{fullpage}
\usepackage{natbib}
\usepackage{setspace}
\usepackage{enumerate}
\usepackage{array}
\usepackage[small]{caption}

\newtheorem{thm}{Theorem}

\newtheorem{prop}[thm]{Proposition}
\newtheorem{cor}[thm]{Corollary}
\numberwithin{equation}{section}
\DeclareMathOperator*{\argmax}{argmax}
\newenvironment{prooftitle}[1]{{\noindent \textsc{Proof #1}}\\}

\begin{document}

\title{Smoothed log-concave maximum likelihood estimation with applications}
\author{Yining Chen and Richard J. Samworth\\ 
  Statistical Laboratory\\
  University of Cambridge\\
  \{y.chen, r.samworth\}@statslab.cam.ac.uk
}

\maketitle

\onehalfspacing

\begin{abstract}
  \noindent
We study the smoothed log-concave maximum likelihood estimator of a probability distribution on $\mathbb{R}^d$.  This is a fully automatic nonparametric density estimator, obtained as a canonical smoothing of the log-concave maximum likelihood estimator.  We demonstrate its attractive features both through an analysis of its theoretical properties and a simulation study.  Moreover, we use our methodology to develop a new test of log-concavity, and show how the estimator can be used as an intermediate stage of more involved procedures, such as constructing a classifier or estimating a functional of the density.  Here again, the use of these procedures can be justified both on theoretical grounds and through its finite sample performance, and we illustrate its use in a breast cancer diagnosis (classification) problem. 


\end{abstract}

\section{Introduction}
\label{Sec:Intro}

Maximum likelihood estimation of shape-constrained densities has received a great deal of interest recently.  The allure is the prospect of obtaining fully automatic nonparametric estimators, with no tuning parameters to choose.  The general idea dates back to \citet{Grenander1956}, who derived the maximum likelihood estimator of a decreasing density on $[0,\infty)$. A characteristic feature of these shape-constrained maximum likelihood estimators is that they are not smooth.  For instance, the Grenander estimator has discontinuities at some of the data points.  The maximum likelihood estimator of a multi-dimensional log-concave density is the exponential of what \citet{CSS2010} call a \emph{tent function}; it may have several ridges. Moreover, in this (and other) examples, the estimator drops discontinuously to zero outside the convex hull of the data.
 
In some applications, the lack of smoothness may not be a drawback in itself.  However, in other circumstances,  a smooth estimate might be preferred, because:
\begin{enumerate}[(a)]
\setlength{\itemsep}{0pt}
\setlength{\parskip}{0pt}
\setlength{\parsep}{0pt}
\item it has a more attractive visual appearance, without ridges or discontinuities that might be difficult to justify to a practitioner;
\item it has the potential to offer substantially improved estimation performance, particularly for small sample sizes, where the convex hull of the data is likely to be rather small;
\item for certain applications, e.g. classification, the maximum likelihood estimator being zero outside the convex hull of the data may present problems; see Section~\ref{Sec:Classification} for further discussion.
\end{enumerate}

For these reasons, we investigate a smoothed version of the $d$-dimensional log-concave maximum likelihood estimator.  The smoothing is achieved by a convolution with a Gaussian density, which preserves the log-concavity shape constraint. To decide how much to smooth, we exploit an interesting property of the log-concave maximum likelihood estimator, which provides a canonical choice of covariance matrix for the Gaussian density, thereby retaining the fully automatic nature of the estimate.  The basic idea, which was introduced by \citet{DuembgenRufibach2009,DuembgenRufibach2011} for the case $d=1$ and touched upon in \citet{CSS2010}, is described in greater detail in Section~\ref{Sec:MDP}.

\begin{figure}[ht!]
  \centering
  $\begin{array}{c c} 
\includegraphics[scale=0.32]{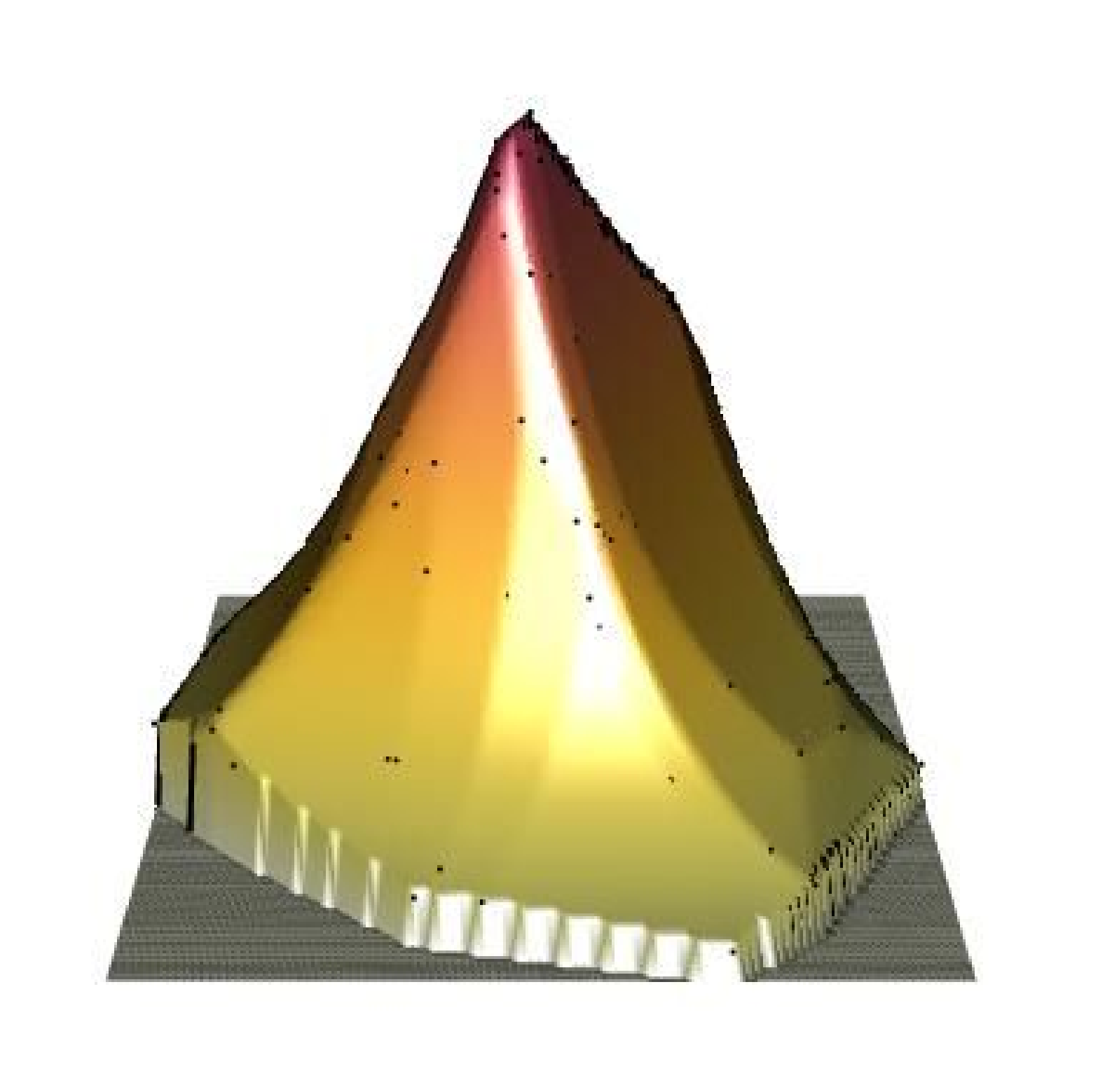} &
\includegraphics[scale=0.32]{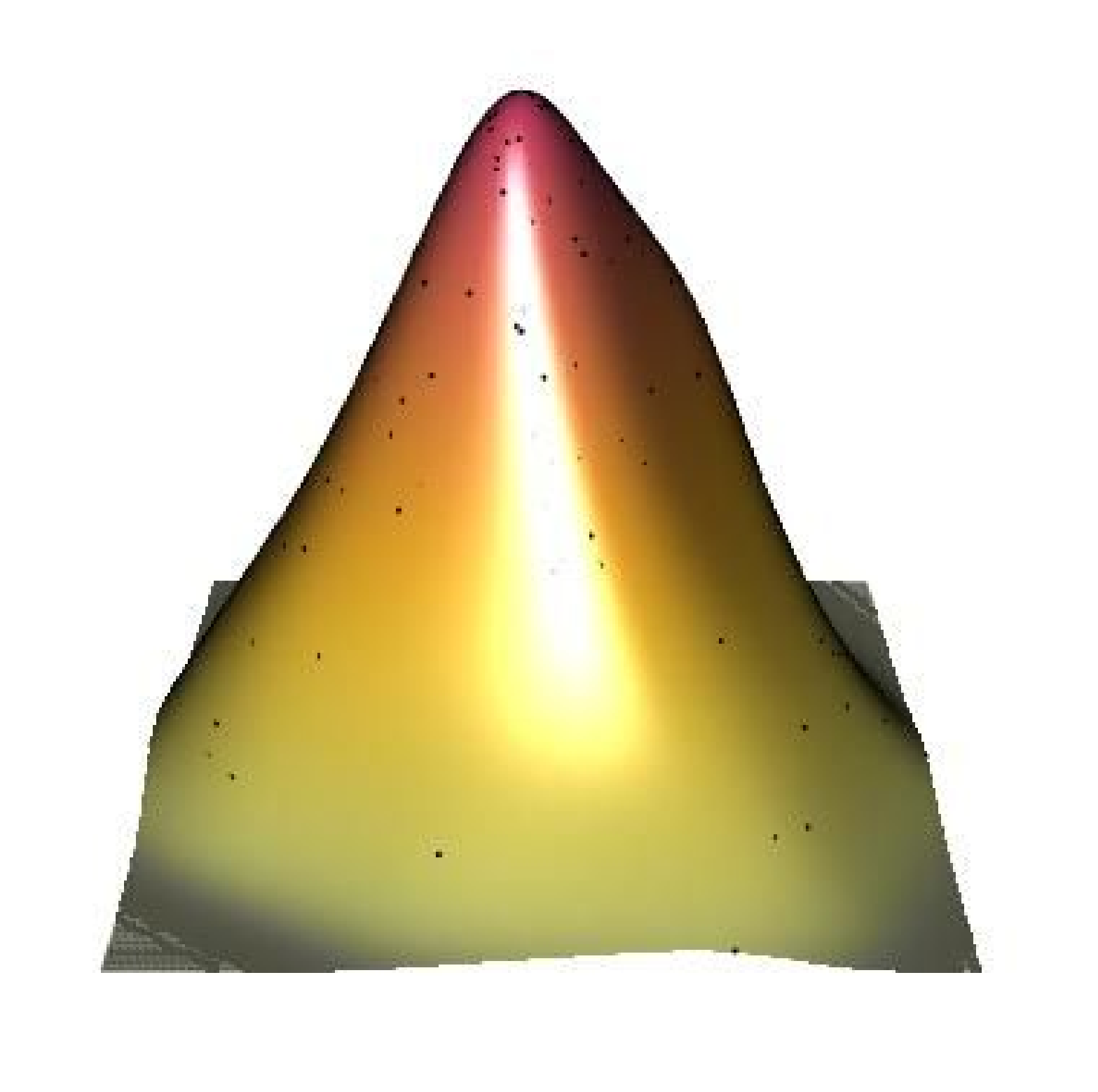} \\
\mathrm{(a)} & \mathrm{(b)} 
  \end{array}$
  \caption{Density estimates based on $n=200$ observations, plotted as dots, from a standard bivariate normal distribution: (a) log-concave maximum likelihood estimator; (b) smoothed log-concave maximum likelihood estimator.}
  \label{fig:lcdsmlcd}
\end{figure}

The challenge of computing the estimator, which involves a $d$-dimensional convolution integral, is taken up in Section~\ref{Sec:Computation}; see Figure~\ref{fig:lcdsmlcd} for an illustration of the estimates obtained.  The theoretical properties of the smoothed log-concave estimator are studied in Section~\ref{Sec:Theory}.  Our framework handles both cases where the log-concavity assumption holds and where it is violated.  
In Section~\ref{Sec:Projections}, we present new results on the infinite-dimensional projection from a probability distribution on $\mathbb{R}^d$ to its closest log-concave approximation; these give further insight into the misspecified setting.  A simulation study follows in Section~\ref{Sec:FSP}, confirming the excellent finite-sample performance.

In Section~\ref{Sec:Test}, we introduce a new hypothesis test of log-concavity of multivariate distributions based on our choice of covariance matrix for the Gaussian density.  This test is consistent, easy to implement, and has much improved finite-sample performance compared to existing methods.   Section~\ref{Sec:Applications} is devoted to applications of the smoothed log-concave maximum likelihood estimator to classification and other functional estimation problems.  We provide theory, under both correct and incorrect model specification, for the performance of the resulting procedures in these cases.  The classification methodology is applied to the Wisconsin breast cancer data set, where the aim is to aid the diagnosis of future potential breast cancer instances.  All proofs are deferred to the Appendix.


Theoretical properties of the unsmoothed log-concave maximum likelihood estimator have been studied in \citet{Walther2002}, \citet{PWM2007}, \citet{BRW2009} and \citet{DuembgenRufibach2009} for the case $d=1$, and \citet{CuleSamworth2010}, \citet{SchuhmacherDuembgen2010} and \citet{DSS2011} for the multivariate case. Further properties of log-concave distributions are discussed in \citet{SHD2011}, and \citet{Walther2009} provides an overview of the field. Other methods for enforcing various shape constraints have been studied in \citet{BraunHall2001}, \citet{GJW2001}, \citet{BalabdaouiWellner2007}, \citet{BalabdaouiWellner2010}, \citet{PavlidesWellner2012}  and \citet{CDH2011}.

\section{The smoothed log-concave maximum likelihood estimator}
\label{Sec:SmoothedLC}

\subsection{Definition and basic properties}
\label{Sec:MDP}

Let $\mathcal{P}$ denote the set of all probability distributions $P$ on $\mathbb{R}^d$ such that $P(H) < 1$ for all hyperplanes $H$.  In this section, we assume that $X_1,X_2,\ldots$ are independent random vectors in $\mathbb{R}^d$ with distribution $P_0 \in \mathcal{P}$.  In that case, for sufficiently large $n$ the convex hull of the data, denoted $C_n = \mathrm{conv}(X_1,\ldots,X_n)$, is $d$-dimensional with probability 1.  It is then known that there exists a unique log-concave density $\hat{f}_n$ that maximises the likelihood function
\[
L(f) = \prod_{i=1}^n f(X_i)
\]
over all log-concave densities $f$.  The estimator $\hat{f}_n$ is supported on $C_n$, and $\log \hat{f}_n$ is piecewise affine on this set.  More precisely, there exists an index set $J$ consisting of $(d+1)$-tuples $j = (j_0,\ldots,j_d)$ of distinct indices in $\{1,\ldots,n\}$, such that $C_n$ can be triangulated into simplices $C_{n,j} = \mathrm{conv}(X_{j_0},\ldots,X_{j_d})$ in such a way that
\[
\log \hat{f}_n(x) = \left\{ \begin{array}{ll} b_j^T x - \beta_j & \mbox{if $x \in C_{n,j}$,} \\
-\infty & \mbox{otherwise,} \end{array} \right.
\]
for some vectors $\{b_j:j \in J\}$ in $\mathbb{R}^d$ and real numbers $\{\beta_j: j \in J\}$.  
Such a function was called a \emph{tent function} in \citet{CSS2010} because when $d=2$ one can think of associating a `tent pole' with each observation, extending vertically out of the plane.  For certain tent pole heights, the graph of $\log \hat{f}_n$ is then the roof of a taut tent stretched over the tent poles.

Despite the attractive asymptotic properties of $\hat{f}_n$ derived in the papers cited in the introduction, the simulation results in \citet{CSS2010} and \citet{Chen2010} indicate that the finite-sample performance is only strong relative to competitors (e.g. kernel-based methods) for moderate or large sample sizes (say $n > 500$).  It appears that for smaller values of $n$, the convex hull of the data is typically not large enough for good performance.  

The idea for fully automatic smoothing of the log-concave maximum likelihood estimator comes from the following observation: Remark~2.3 of \citet{DSS2011} (see also Corollary~2.3 of \citet{DuembgenRufibach2009}) shows that while the log-concave maximum likelihood estimator is a good estimator of the first moment of $P_0$, it underestimates the covariance matrix.  More precisely, we have that 
\[
\int_{\mathbb{R}^d} x \hat{f}_n(x) \, dx = \frac{1}{n}\sum_{i=1}^n X_i \equiv \bar{X},
\]
say.  On the other hand, however, 
\begin{align}
\label{Eq:Sigmas}
\tilde{\Sigma} \equiv \int_{\mathbb{R}^d} (x - \bar{X})(x - \bar{X})^T \hat{f}_n(x) \, dx &\leq \frac{1}{n}\sum_{i=1}^n (X_i - \bar{X})(X_i - \bar{X})^T \nonumber \\
&< \frac{1}{n-1}\sum_{i=1}^n (X_i - \bar{X})(X_i - \bar{X})^T  \equiv \hat{\Sigma}.
\end{align}
Here, $A \leq B$ and $A <B$ mean the matrix $B - A$ is non-negative definite and positive definite respectively.

This allows us to define our modified estimator, which we call the \emph{smoothed log-concave maximum likelihood estimator} and denote $\tilde{f}_n$.  It is given by
\begin{equation}
\label{Eq:fntilde}
\tilde{f}_n = \hat{f}_n \ast \phi_{d,\hat{A}},
\end{equation}
where $\phi_{d,\hat{A}}$ is the $d$-variate normal density with zero mean and covariance matrix $\hat{A} = \hat{\Sigma} - \tilde{\Sigma}$.  Note that the level of smoothing is automatically determined through the matrix $\hat{A}$.  

The basic properties of $\tilde{f}_n$ are summarised in the proposition below.
\begin{prop}
\label{Prop:Basic}
Let $P_0 \in \mathcal{P}$, and let $\tilde{f}_n$ denote the smoothed log-concave maximum likelihood estimator $\tilde{f}_n$ based on independent observations $X_1,\ldots,X_n$ having distribution $P_0$.  Then
\renewcommand{\labelenumi}{(\theenumi)}
\begin{description}
\item[(a)] $\tilde{f}_n$ is log-concave;
\item[(b)] the support of $\tilde{f}_n$ is $\mathbb{R}^d$;
\item[(c)] $\tilde{f}_n$ is a real analytic function on $\mathbb{R}^d$ (in particular, it is infinitely differentiable);
\item[(d)] the mean and covariance matrix corresponding to $\tilde{f}_n$ agree with the sample mean and sample covariance matrix: $\int_{\mathbb{R}^d} x \tilde{f}_n(x) \, dx = \bar{X}$ and $\int_{\mathbb{R}^d} (x - \bar{X})(x - \bar{X})^T \tilde{f}_n(x) \, dx = \hat{\Sigma}$.
\end{description}
\end{prop} 

\subsection{Computational issues}
\label{Sec:Computation}

The aim of this section is to describe algorithms for computing the smoothed log-concave maximum likelihood estimator $\tilde{f}_n$.  As a preliminary step, we need to compute the covariance matrix $\hat{A}$ of the multivariate normal distribution used in the convolution~(\ref{Eq:fntilde}).

\subsubsection{Computation of the covariance matrix $\hat{A}$}
\label{Sec:CovMatrix}

Recall that $\hat{A} = \hat{\Sigma} - \tilde{\Sigma}$, where $\hat{\Sigma}$ is the sample covariance matrix, and 
\begin{equation}
\label{Eq:Sigmatilde}
\tilde{\Sigma}= \int_{\mathbb{R}^d} xx^T \, \hat{f}_n(x) \, dx - \bar{X}\bar{X}^T = \sum_{j \in J} \int_{C_{n,j}} xx^T \exp(b_j^T x - \beta_j) \, dx - \bar{X}\bar{X}^T.
\end{equation}
We make an affine transformation of each of the regions of integration onto the unit simplex.  Recall that $C_{n,j} = \mathrm{conv}(X_{j_0},\ldots,X_{j_d})$, set $D_j= \det[X_{j_1}-X_{j_0},X_{j_2}-X_{j_0},\ldots,X_{j_d}-X_{j_0}]$, and let $U_d=\{u = (u_1,\ldots,u_d) \in [0,\infty)^d:\sum_{l=1}^d u_l \leq 1\}$ be the unit simplex in $\mathbb{R}^d$.  Following \citet{CuleDuembgen2008}, we further define the auxiliary functions $J_{\tilde{d}}:\mathbb{R}^{\tilde{d}+1} \rightarrow \mathbb{R}$ by
\[
J_{\tilde{d}}(y_0,y_1,\ldots,y_{\tilde{d}}) = \int_{U_{\tilde{d}}} \exp\biggl(\sum_{l=0}^{\tilde{d}} u_l y_l\biggr) \, du_1 \ldots du_{\tilde{d}},
\]
where $u_0 = 1 - \sum_{l=1}^{\tilde{d}} u_l$.  Then, writing $y_{j_l} = \log \hat{f}_n(X_{j_l})$, we have
\begin{align*}
\sum_{j \in J} \int_{C_{n,j}} xx^T \exp(b_j^T x - \beta_j) \, dx &= \sum_{j \in J} |D_j| \int_{U_d} \biggl(\sum_{l=0}^d u_l X_{j_l}\biggr)\biggl(\sum_{l=0}^d u_l X_{j_l}\biggr)^T e^{\sum_{l=0}^d u_l y_{j_l}} \, du \\
&= \sum_{j \in J} |D_j| \biggl\{\sum_{l=0}^d \sum_{l'=0}^d X_{j_l}X_{j_{l'}}^T \frac{\partial^2 J_d(y_{j_0},y_{j_1},\ldots,y_{j_d})}{\partial y_{j_l} \partial y_{j_{l'}}}\biggr\} \\
&= \sum_{j \in J} |D_j| \biggl\{ \sum_{l=0}^d \sum_{l'=0}^d X_{j_l}X_{j_{l'}}^T J_{d+2}(y_{j_0},y_{j_1},\ldots,y_{j_d},y_{j_l},y_{j_{l'}}) \\
&\hspace{3cm}+ \sum_{l=0}^d X_{j_l}X_{j_l}^T J_{d+2}(y_{j_0},y_{j_1},\ldots,y_{j_d},y_{j_l},y_{j_l}) \biggr\}.
\end{align*}
We have applied the basic results of \citet{CuleDuembgen2008} in the last step.  An exact expression for $J_{d+2}(\cdot)$ is given in Appendix~B.1 of \citet{CSS2010} when its arguments are non-zero and distinct.  The Taylor approximation of \citet{CuleDuembgen2008} can be used when some of the arguments are small or have similar (or equal) values.

\subsubsection{Computation of the smoothed log-concave maximum likelihood estimator}

We have
\[
\tilde{f}_n(x_0)=\sum_{j \in J} \int_{C_{n,j}}e^{b_j^Tx-\beta_j} \frac{1}{(2\pi)^{d/2}(\det \hat{A})^{1/2}}\, e^{-\frac{1}{2}(x_0-x)^T\hat{A}^{-1}(x_0-x)} \, dx.
\]
By making an affine transformation of each ${C_{n,j}}$ onto the unit simplex as in Section~\ref{Sec:CovMatrix}, we reduce the problem to integrating the exponential of a quadratic polynomial over the unit simplex. In general, this has no explicit solution, so it has to be evaluated numerically.

\citet{Stroud1971} gives a brief introduction to the problem of evaluating integrals over the unit simplex, while \citet{GrundmannMoeller1978} proposed a combinatorial method.  We apply their method, first noting that by integrating out one variable, the dimensionality of the integral can be reduced by one.  To see this, consider any $d \times d$ positive definite, symmetric matrix $A \equiv [a_{ll'}]$, any vector $B = (b_1,\ldots,b_d)^T \in \mathbb{R}^d$ and any constant $c \in \mathbb{R}$.  Writing $\Phi(\cdot)$ for the standard normal distribution function, $u = (u_1,\ldots,u_d)^T$ and $u_0 = 1 - \sum_{l=1}^{d-1} u_l$ we have  
\begin{align}
& \int_0^1 \!\int_0^{1-u_1} \! \cdots \! \int_0^{1-u_1- \cdots -u_{d-1}} \! e^{-u^TAu+B^Tu+c} \, du_d \ldots du_2 du_1 \notag\\ 
&= \int_0^1 \! \int_0^{1-u_1} \! \cdots \! \int_0^{1-u_1- \cdots -u_{d-1}} \! e^{-a'u_d^2+b'u_d+c'} du_d \ldots du_2 du_1 \notag\\
&= \int_0^1 \! \int_0^{1-u_1} \! \cdots \! \int_0^{1-u_1- \cdots -u_{d-2}} \! e^{c'+\frac{b'^2}{4a'}}\sqrt{\frac{\pi}{a'}}\biggl\{\Phi\Bigl(u_0\sqrt{2a'}\!-\!\frac{b'}{\sqrt{2a'}}\Bigr) - \Phi\Bigl(\frac{-b'}{\sqrt{2a'}}\Bigr)\biggr\} \, du_{d-1} \ldots du_2 du_1. \label{Eq:Integ1D}
\end{align} 
Here, $a'$, $b'$ and $c'$ are defined by
\[
a' = a_{dd}, \quad b' = b_d+2\sum_{l=1}^{d-1}a_{dl}u_l \quad \text{and} \quad c' = u_{-d}^T[a_{ll'}]_{1 \leq l,l' \leq d-1}u_{-d}+\sum_{l=1}^{d-1} b_lu_l + c,
\]
where $u_{-d} = (u_1,\ldots,u_{d-1})^T$.  It follows that we can use the combinatorial method to integrate over the $(d-1)$-dimensional unit simplex.  Some special cases include:
\begin{enumerate}[(a)]
\setlength{\itemsep}{0pt}
\setlength{\parskip}{0pt}
\setlength{\parsep}{0pt}
\item $d = 1$. In this case, (\ref{Eq:Integ1D}) is a simple function of $\Phi(\cdot)$, and the smoothed log-concave maximum likelihood estimator can be computed straightforwardly. This method is implemented in the \texttt{R} package \texttt{logcondens} \citep{RufibachDuembgen2006,DuembgenRufibach2011}.
\item $d = 2$. In this case, (\ref{Eq:Integ1D}) is an integral over $[0,1]$, and other standard numerical integration methods such as the Gaussian quadrature rule, can be applied.  
\end{enumerate}

The combinatorial method and its variations are implemented in the latest version of the \texttt{R} package \texttt{LogConcDEAD} \citep{CGSC2007,CGS2009}.  We found this method to be numerically stable even with several thousand observations, when $\det \hat{A}$ may be rather small (note that in such cases, $a'$ in~(\ref{Eq:Integ1D}) will typically not be close to zero).  However, we briefly present below two other ways of computing $\tilde{f}_n(x_0)$; while slower in most cases, they do not require the inversion of $\hat{A}$, so can be used even when $\det \hat{A}$ is very small.
\begin{enumerate}[(a)]
\setlength{\itemsep}{0pt}
\setlength{\parskip}{0pt}
\setlength{\parsep}{0pt}
\item \textbf{Monte Carlo method}. 
\begin{enumerate}[(1)]
\setlength{\itemsep}{0pt}
\setlength{\parskip}{0pt}
\setlength{\parsep}{0pt}
\item Conditional on $X_1,\ldots,X_n$, generate independent random vectors $X_1^*, \ldots, X_B^*$ from the $N_d(x_0,\hat{A})$ distribution.
\item Approximate $\tilde{f}_n(x_0)$ by $\frac{1}{B} \sum_{b=1}^B \hat{f}_n(X_b^*)$.
\end{enumerate}
The validity of this approximation follows from the strong law of large numbers, applied conditional on $X_1,\ldots,X_n$.
\item \textbf{Fourier transform}. We can take advantage of the convolution property of the Fourier transform $\mathcal{F}$ as follows.  First note that
\[
\mathcal{F}(\hat{f}_n)(\xi) = \int_{\mathbb{R}^d} \hat{f}_n(x) e^{-i \xi^T x} dx = \sum_{j \in J} \int_{C_{n,j}}e^{(b_j- i \xi)^T x-\beta_j} dx,
\] 
which can be evaluated by extending the auxiliary functions $J_d$ to the complex plane.  Since $\mathcal{F}(\tilde{f}_n)(\xi) = e^{-i x_0^T \xi - \xi^T\hat{A}\xi/2}\mathcal{F}(\hat{f}_n)(\xi)$, we can invert $\mathcal{F}(\tilde{f}_n)$ on a fine grid using the fast Fourier transform.
\end{enumerate}

\subsubsection{Sampling from the fitted density estimate}
\label{Sec:Sampling}

Since $\tilde{f}_n$ is the convolution of $\hat{f}_n$ and a multivariate normal density, conditional on $X_1,\ldots,X_n$, it is straightforward to draw an observation $X^{**}$ from ${\tilde{f}_n}$ as follows:
\begin{enumerate}[(a)]
\setlength{\itemsep}{0pt}
\setlength{\parskip}{0pt}
\setlength{\parsep}{0pt}
\item Draw $X^*$ from ${\hat{f}_n}$ using the algorithm described in Appendix~B.3 of \citet{CSS2010} or the algorithm of \citet{GopalCasella2010}.
\item Draw $u \sim N_d(0,\hat{A})$, independent of $X^*$.
\item Return $X^{**} = X^* + u$.
\end{enumerate}

\subsection{Theoretical performance}
\label{Sec:Theory}

It is convenient to define, for $r=1,2$, the classes of probability distributions on $\mathbb{R}^d$ given by 
\[
\mathcal{P}_r = \biggl\{P \in \mathcal{P}:\int_{\mathbb{R}^d} \|x\|^r \, dP(x) < \infty\biggr\}.
\]
The condition $P_0 \in \mathcal{P}_1$ is necessary and sufficient for the existence of a unique upper semi-continuous log-concave density $f^*$ that maximises $\int \log f \, dP_0$ over all log-concave densities $f$ \citep[Theorem~2.2]{DSS2011}.  In fact, if $P_0$ has a density $f_0$, and provided that $\int f_0 \log f_0 < \infty$ (which is certainly the case if $f_0$ is bounded), $f^*$ minimises the Kullback--Leibler divergence $d_{KL}(f,f_0) = \int f_0 \log (f_0/f)$ over all log-concave densities $f$.  In this sense, $f^*$ is the closest log-concave density to $P_0$.  

The density $f^*$ plays an important role in the following theorem, which describes the asymptotic behaviour of the smoothed log-concave maximum likelihood estimator $\tilde{f}_n$.
\begin{thm}
\label{Thm:Asymp}
Suppose that $P_0 \in \mathcal{P}_2$, and write $\mu = \int_{\mathbb{R}^d} x \, dP_0(x)$ and $\Sigma = \int_{\mathbb{R}^d} (x-\mu)(x-\mu)^T \, dP_0(x)$.  Let $f^{**} = f^* \ast N_d(0,A^*)$, where $A^* = \Sigma - \Sigma^*$ with $\Sigma^* = \int_{\mathbb{R}^d}(x-\mu)(x-\mu)^T f^*(x) \, dx$.  Taking $a_0 > 0$ and $b_0 \in \mathbb{R}$ such that $f^{**}(x) \leq e^{-a_0\|x\| + b_0}$, we have for all $a < a_0$ that
\[
\int_{\mathbb{R}^d} e^{a\|x\|}|\tilde{f}_n(x) - f^{**}(x)| \stackrel{a.s.}{\rightarrow} 0
\]
and, if $f^{**}$ is continuous, $\sup_{x \in \mathbb{R}^d} e^{a\|x\|}|\tilde{f}_n(x) - f^{**}(x)| \stackrel{a.s.}{\rightarrow} 0$. 
\end{thm}
The condition that $P_0 \in \mathcal{P}_2$ imposed in Theorem~\ref{Thm:Asymp} ensures the finiteness of $A^*$. We see that in general, $\tilde{f}_n$ converges to a slightly smoothed version of the closest log-concave density to $P_0$.  However, if $P_0$ has a log-concave density $f_0$, then $f_0 = f^* = f^{**}$, so $\tilde{f}_n$ is strongly consistent in these exponentially weighted total variation and supremum norms.  In fact, suppose that $a:\mathbb{R}^d \rightarrow \mathbb{R}$ is a sublinear function, i.e. $a(x+y) \leq a(x) + a(y)$ and $a(rx) = ra(x)$ for all $x,y \in \mathbb{R}^d$ and $r \geq 0$, satisfying $e^{a(x)}f(x) \rightarrow 0$ as $\|x\| \rightarrow \infty$.  It can be shown that under the conditions of Theorem~\ref{Thm:Asymp},
\[
\int_{\mathbb{R}^d} e^{a(x)}|\tilde{f}_n(x) - f^{**}(x)| \stackrel{a.s.}{\rightarrow} 0
\]
\citep{SHD2011}.

Despite being smooth and having full support, it turns out that $\tilde{f}_n$ is rather close to $\hat{f}_n$.  This is quantified in the finite-sample bound below.
\begin{prop}
\label{Prop:Bounds}
\begin{description}
\item If $x \in C_{n,j}$, and $\hat{f}_n(x) = \exp(b_j^T x - \beta_j)$, then
\[
\frac{\tilde{f}_n(x) - \hat{f}_n(x)}{\hat{f}_n(x)} \leq e^{\frac{1}{2}b_j^T \hat{A}b_j} - 1.
\]
Moreover,   
\[
\int_{\mathbb{R}^d} |\tilde{f}_n - \hat{f}_n| \leq 2(e^{\frac{1}{2}\lambda_{\max}} - 1 + \delta_n) 
\]
where $\lambda_{\max} = \max_{j \in J} b_j^T \hat{A} b_j$, and $\delta_n = \int_{C_n^c} \tilde{f}_n$.
\end{description}
\end{prop}

\subsection{Properties of (smoothed) log-concave approximations}
\label{Sec:Projections}

In this subsection, we give new insights into the maps from a probability distribution $P$ to its log-concave approximation $f^*$, and its smoothed version $f^{**}$. Results such as these enhance our understanding of the behaviour of maximum likelihood estimators in non-convex, misspecified models, where existing results are very limited.  Theorem~\ref{Thm:Independent} below shows that log-concave approximations and their smoothed analogues preserve independence of components.  As well as being of use in our simulation studies, this is the key result which underpins a new approach to fitting independent component analysis models using nonparametric maximum likelihood \citep{SamworthYuan2012}.  
\begin{thm}
\label{Thm:Independent}
Suppose that $P \in \mathcal{P}_1$ is a product measure on $\mathbb{R}^d$, so that $P = P_1 \otimes P_2$, say, where $P_1$ and $P_2$ are probability measures on $\mathbb{R}^{d_1}$ and $\mathbb{R}^{d_2}$ respectively, with $d_2 = d -d_1$.   Let $f^*$ denote the log-concave approximation to $P$, and let $f_\ell^*$ denote the log-concave approximation to $P_\ell$, for $\ell=1,2$.  Then, writing $x = (x_1^T,x_2^T)^T$, where $x_1 \in \mathbb{R}^{d_1}$ and $x_2 \in \mathbb{R}^{d_2}$, we have
\[
f^*(x) = f_1^*(x_1)f_2^*(x_2).
\]
Now suppose further that $P \in \mathcal{P}_2$.  Let $f^{**}$ denote the smoothed log-concave approximation to $P$, and let $f_\ell^{**}$ denote the smoothed log-concave approximation to $P_\ell$, for $\ell=1,2$.  Then, for all $x = (x_1^T,x_2^T)^T$,
\[
f^{**}(x) = f_1^{**}(x_1)f_2^{**}(x_2).
\]
\end{thm}

Our next theorem characterises the log-concavity constraint through the trace of the non-negative definite matrix $A^*$ defined in Theorem~\ref{Thm:Asymp}.  
\begin{thm}
\label{Thm:Covariance}
Suppose that $P \in \mathcal{P}_1$.  Then $\operatorname{tr}(A^*) = 0$ if and only if $P$ has a log-concave density. 
\end{thm}
The `if' part of this statement is well-known, but the `only if' part is new.  The two parts together motivate our testing procedure for log-concavity, which is developed in Section~\ref{Sec:Test}.

In most cases, it is very difficult to find explicitly the log-concave approximation $f^*$ to a given distribution $P \in \mathcal{P}_1$.  Our final result of this section is straightforward to prove, but is of interest because it shows that some log-concave densities can have a large `domain of attraction'.
\begin{prop}
\label{Prop:Convex}
Let $f^*$ be an upper semi-continuous, log-concave density on $\mathbb{R}^d$.  Then the class of distributions $P \in \mathcal{P}_1$ with log-concave approximation $f^*$ is convex.
\end{prop}
For instance, if $f(x;\alpha,\sigma) = \frac{\alpha \sigma^\alpha}{2(|x|+\sigma)^{\alpha+1}}$ is a symmetrised Pareto density with $\alpha > 1$ and $\sigma > 0$, then it can be shown that its log-concave projection is $f^*(x;\alpha,\sigma) = \frac{\alpha-1}{2\sigma}\exp\{-(\alpha - 1)|x|/\sigma\}$.  Thus the class of distributions with whose log-concave projection is the standard Laplace density is infinite-dimensional.

\subsection{Finite sample performance}
\label{Sec:FSP}

\begin{figure}[hp!]
  \centering
\includegraphics[scale=0.70]{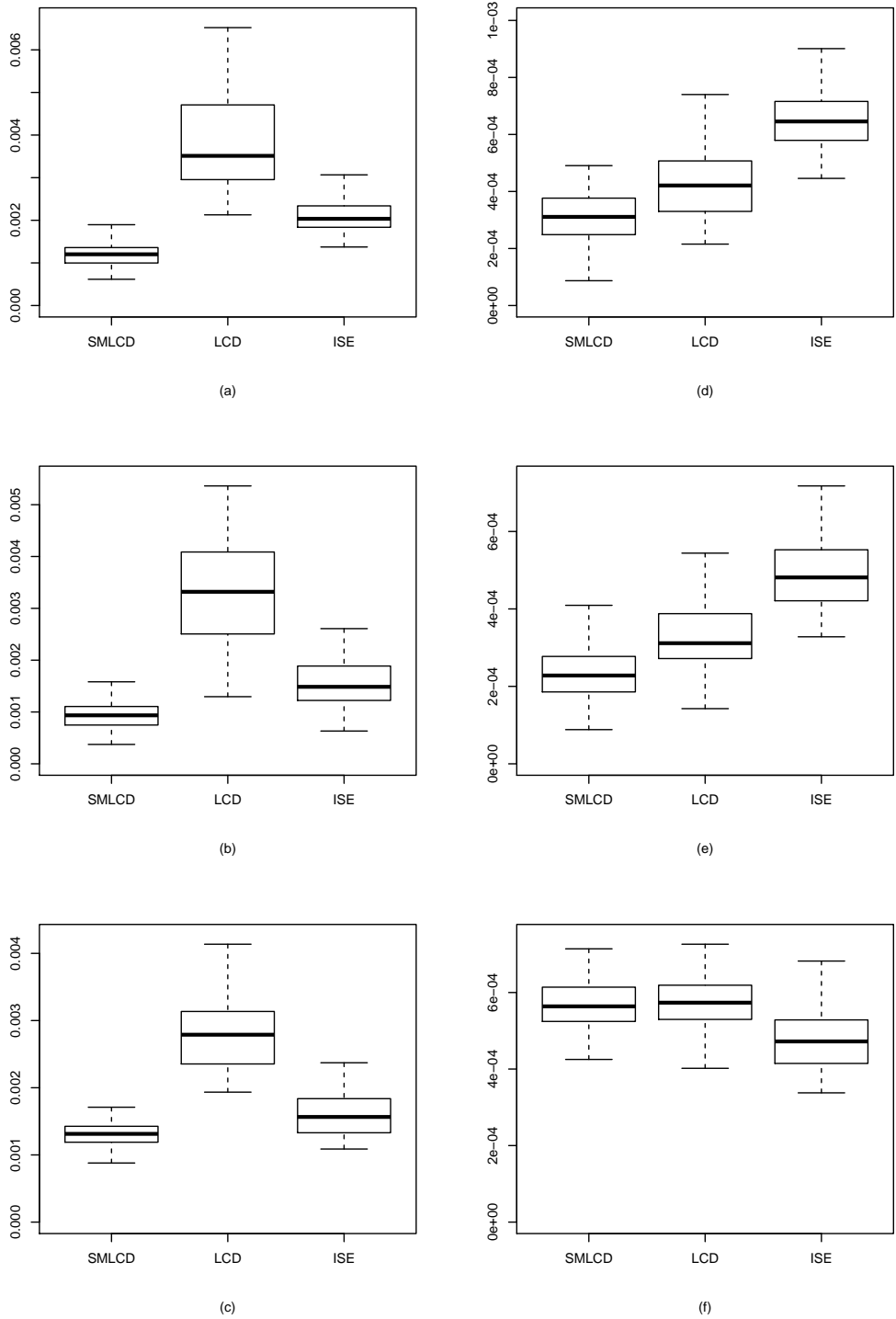}
  \caption{Boxplots of ISEs for $d=3$ with the Gaussian location mixture true density for the smoothed log-concave maximum likelihood estimator SMLCD, log-concave maximum likelihood estimator LCD and kernel density estimator with the `oracle' optimal ISE bandwidth: (a) $n=100$, $\|\mu\|=1$; (b) $n=100$, $\|\mu\|=2$; (c) $n=100$, $\|\mu\|=3$; (d) $n=1000$, $\|\mu\|=1$; (e) $n=1000$, $\|\mu\|=2$; (f) $n=1000$, $\|\mu\|=3$.}
  \label{Fig:box3d}
\end{figure}

Our simulation study considered the normal location mixture density $f(\cdot) = 0.4\phi_d(\cdot)+0.6\phi_d(\cdot-\mu)$ for $\|\mu\|=$ 1, 2 and 3, where $\phi_d = \phi_{d,I}$.  This mixture density is log-concave if and only if  $\|\mu\| \leq 2$.  For each density, for $d=2$ and $d=3$, and for sample sizes $n=100$ and $n=1000$, we computed the Integrated Squared Error (ISE) of the smoothed log-concave maximum likelihood estimator for each of 50 replications. We also computed the ISE of the log-concave maximum likelihood estimator and that of a kernel density estimator with a Gaussian kernel and the optimal ISE bandwidth for each individual data set, which would be unknown in practice. The boxplots of the ISEs for the different methods are given in Figure~\ref{Fig:box3d} for $d=3$. The analogous plots for the case $d=2$ can be found in \citet{ChenSamworth2011}.


We see that when the true density is log-concave, the smoothed log-concave estimator offers substantial ISE improvements over its unsmoothed analogue for both sample sizes, particularly at the smaller sample size $n=100$. It also outperforms by a considerable margin the kernel density estimator with the optimal ISE bandwidth.  When the log-concavity assumption is violated, the smoothed log-concave estimator is still competitive with the optimal-ISE kernel estimator at the smaller sample size $n=100$, and also improves on its unsmoothed analogue.  However, at the larger sample size $n=1000$, the bias caused by the fact that $\int_{\mathbb{R}^d} (f^* - f)^2 > 0$ dominates the contribution from the variance of the estimator, and the kernel estimator is an improvement. These results confirm that the smoothed log-concave estimator has excellent performance when the true density is log-concave, and remains competitive in situations where the log-concavity assumption is violated, provided that the modelling bias caused by this misspecification is not too large relative to the sampling variability of the estimator.

\section{A new test of log-concavity}
\label{Sec:Test}

Several tests of log-concavity have been proposed in the literature.  \citet{An1995} and \citet{Walther2002} discuss various tests for univariate data, while \citet{CSS2010} presented two tests of log-concavity for multivariate data.  \citet{Hazelton2011} proposed another multivariate test based on kernel density estimates which had improved finite-sample performance on his simulated examples.  However, none of these multivariate tests has theoretical support.

Suppose $X_1,\ldots,X_n \stackrel{iid}{\sim} P_0 \in \mathcal{P}_1$, and we seek a size $\alpha \in (0,1)$ test of $H_0: P_0$ has a log-concave density against $H_1: P_0$ does not have a log-concave density.  Motivated by Theorem~\ref{Thm:Covariance}, we propose the following procedure:
\begin{enumerate}[(a)]
\setlength{\itemsep}{0pt}
\setlength{\parskip}{0pt}
\setlength{\parsep}{0pt}
\item Compute the log-concave maximum likelihood density estimate $\hat{f}_n$.
\item Compute the test statistic $\operatorname{tr}(\hat{A})$, where $\hat{A} = \hat{\Sigma} - \tilde{\Sigma}$, as in~(\ref{Eq:Sigmas}).
\item Generate a reference distribution as follows: for $b=1,\ldots,B$, draw conditionally independent samples $X_{1b}^*,\ldots,X_{nb}^*$ from $\hat{f}_n$.  For each bootstrap sample, first compute the log-concave maximum likelihood estimator $\hat{f}_{nb}$.  Then compute $\operatorname{tr}(\hat{A}_{nb})$, where
\[
\hat{A}_{nb} \equiv \hat{\Sigma}_b - \tilde{\Sigma}_b \equiv \frac{1}{n-1}\sum_{i=1}^n (X_{ib}^* - \bar{X}_b^*)(X_{ib}^* - \bar{X}_b^*)^T - \int_{\mathbb{R}^d} (x - \bar{X}_b^*)(x - \bar{X}_b^*)^T \hat{f}_{nb}(x) \, dx,
\]
and $\bar{X}_b^* = n^{-1}\sum_{i=1}^n X_{ib}^*$.
\item Reject $H_0$ if $(B+1)^{-1}\sum_{b=1}^{B+1} \mathbbm{1}_{\{\operatorname{tr}(\hat{A}) > \operatorname{tr}(\hat{A}_{nb})\}} > 1 - \alpha$.
\end{enumerate}
We call this procedure a \emph{trace} test. It is justified by the following result: 
\begin{thm}
\label{Thm:Test}
Suppose that $P_0 \in \mathcal{P}_1$. The trace test is consistent: that is, if $P_0$ is not log-concave, then for each $B \in \mathbb{N}$, the power of the test converges to one as $n \rightarrow \infty$.
\end{thm}
We remark that if $P_0 \in \mathcal{P}_2$, one can also draw bootstrap samples from $\tilde{f}_n$ instead of $\hat{f}_n$ in Step (c).  To illustrate the performance of the test, we ran two small simulation studies. In the first study, we simulated from the bivariate mixture of normal distributions density $f(x) = \frac{1}{2}\phi_{2,I}(x) + \frac{1}{2}\phi_{2,I}(x-\mu)$, with $\|\mu\| = 0,2,4$ (which we recall is log-concave if and only if $\|\mu\| \leq 2$). For each simulation setup, we performed 200 hypothesis tests with $B=99$.  The proportion of times that the null hypothesis was rejected in a size $\alpha=0.05$ test is reported in Table~\ref{tab:sim1}. For comparison, we also report the results from the critical bandwidth test proposed by \citet{Hazelton2011}.  The permutation test studied by \citet{CSS2010} did not perform as well as the critical bandwidth test, so we omitted its results here.
\begin{table}[ht]
    \centering
    \begin{tabular}{ | c c c c c |}
    \hline 
    $n$ & Method & $\|\mu\|=0$ & $\|\mu\|=2$ & $\|\mu\|=4$ \\ \hline
    200 & critical bandwidth	& 0.065 & 0.015 & 0.985\\
	& trace			& 0.045 & 0.045 & 1.000\\
    500 & critical bandwidth	& 0.045 & 0.005 & 1.000\\
	& trace			& 0.045 & 0.055 & 1.000\\
    \hline
    \end{tabular}
    \caption{Proportion of times out of 200 repetitions that the null hypothesis was rejected with $\alpha = 0.05$.}
    \label{tab:sim1}
\end{table}
For the second study, we replicate the settings considered in \citet{Hazelton2011}, where four different types of bivariate densities of independent components were chosen.  The marginal distributions were: 
\begin{enumerate}[(a)]
\setlength{\itemsep}{0pt}
\setlength{\parskip}{0pt}
\setlength{\parsep}{0pt}
\item A $\frac{1}{2}N(0,1/4) + \frac{1}{2}N(0,4)$ distribution and a $\frac{1}{2}N(0,1/4) + \frac{1}{2}N(2,4)$ distribution;  \item A $t_4$ distribution in both cases;
\item A $\frac{1}{2}N(0,1/4) + \frac{1}{2}N(2,4)$ and a $t_4$ distribution; 
\item A $\frac{1}{2}N(0,1/4) + \frac{1}{2}N(2,5)$ density, and a $\Gamma(2,1)$ distribution. 
\end{enumerate}
Note that all of these densities are unimodal but not log-concave.  The corresponding estimates of the power of the tests are presented in Table~\ref{tab:sim2}. 
\begin{table}[ht]
    \centering
    \begin{tabular}{ | c c c c c c |}
    \hline 
    $n$ & Method & \multicolumn{4}{c|}{Cases} \\ \cline{3-6}
        & 			& (a)   & (b)   & (c)   & (d) \\ \cline{1-6}
    200 & critical bandwidth 	& 0.520 & 0.195 & 0.395 & 0.295\\
	& trace		 	& 1.000 & 0.960 & 1.000 & 1.000\\
    500 & critical bandwidth 	& 0.760 & 0.340 & 0.710 & 0.505\\
	& trace		 	& 1.000 & 1.000 & 1.000 & 1.000\\
    \hline
    \end{tabular}
    \caption{Proportion of times out of 200 repetitions that the null hypothesis was rejected with $\alpha = 0.05$.}
    \label{tab:sim2}
\end{table}
The first study confirms that the trace test controls the Type I error satisfactorily (and appears to be less conservative than the critical bandwidth test when $\|\mu\| = 2$).  The results of the second study, though, are quite striking, and suggest that our new test for log-concavity has considerably improved finite-sample power compared to the critical bandwidth test.  \citet{Hazelton2011} noted that the critical bandwidth test can have reduced power due to the boundary bias of the kernel estimators and is quite sensitive to the outliers (in fact, one also needs to pick a compact region containing the majority of the data, and this choice is somewhat arbitrary).  Our test avoids these issues and performs well even in the presence of outliers or when the true density has bounded support.

\section{Other applications}
\label{Sec:Applications}

\subsection{Classification problems}
\label{Sec:Classification}

Changing notation slightly from the previous section, we now assume that $(X,Y)$, $(X_1,Y_1)$, $\ldots$, $(X_n,Y_n)$ are independent and identically distributed pairs taking values in $\mathbb{R}^d \times \{1,\ldots,K\}$.  Let $\mathbb{P}(Y = k) = \pi_k$ for $k=1,\ldots,K$, and suppose that conditional on $Y=k$, the random vector $X$ has distribution $P_k$.  

A \emph{classifier} is a measurable function $C:\mathbb{R}^d \rightarrow \{1,\ldots,K\}$, with the interpretation that the classifier assigns the point $x \in \mathbb{R}^d$ to class $C(x)$.  The \emph{misclassification error rate}, or \emph{risk}, of $C$ is 
\[
\mathrm{Risk}(C) = \mathbb{P}\{C(X) \neq Y\}.
\]
In the case where each distribution $P_k$ has a density $f_k$, the classifier that minimises the risk is the \emph{Bayes classifier} $C^{\mathrm{Bayes}}$, given by
\[
C^{\mathrm{Bayes}}(x) = \argmax_{k \in \{1,\ldots,K\}} \pi_k f_k(x).
\]
(For all classifiers defined by an $\argmax$ as above, we will for the sake of definiteness split ties by taking the smallest element of the $\argmax$.)  We will also be interested in the \emph{log-concave Bayes classifier} and \emph{smoothed log-concave Bayes classifier}, defined respectively by
\[
C^{\mathrm{LCBayes}}(x) = \argmax_{k \in \{1,\ldots,K\}} \pi_k f_k^*(x) \quad \text{and} \quad C^{\mathrm{SLCBayes}}(x) = \argmax_{k \in \{1,\ldots,K\}} \pi_k f_k^{**}(x).
\]
Here, $f_k^*$ and $f_k^{**}$ are the log-concave approximation to $P_k$ and its smoothed analogue, defined in Theorem~\ref{Thm:Asymp}.  In particular, both classifier coincide with the Bayes classifier when $\{P_k:k = 1,\ldots,K\}$ have log-concave densities.  Empirical analogues of these theoretical classifiers are given by
\[
\hat{C}_n^{\mathrm{LC}}(x) = \argmax_{k \in \{1,\ldots,K\}} N_k \hat{f}_{n,k}(x) \quad \text{and} \quad \hat{C}_n^{\mathrm{SLC}}(x) = \argmax_{k \in \{1,\ldots,K\}} N_k \tilde{f}_{n,k}(x).
\]
Here, $N_k = \sum_{i=1}^n \mathbbm{1}_{\{Y_i = k\}}$ is the number of observations from the $k$th class, and $\hat{f}_{n,k}$ and $\tilde{f}_{n,k}$ are respectively the log-concave maximum likelihood estimator of $f_k$ and its smoothed analogue, based on $\{X_i: Y_i = k\}$.  

The theorem below describes the asymptotic behaviour of these classifiers.  It reveals that the risk of $\hat{C}_n^{\mathrm{LC}}$ and $\hat{C}_n^{\mathrm{SLC}}$ converges not (in general) to the Bayes risk, but instead to the risk of $C^{\mathrm{LCBayes}}$ and $C^{\mathrm{SLCBayes}}$ respectively.  This is a similar situation to that encountered when a parametric classifier such as linear or quadratic discriminant analysis is used, but the relevant parametric modelling assumptions fail to hold.  It suggests that the classifiers $\hat{C}_n^{\mathrm{LC}}$ and $\hat{C}_n^{\mathrm{SLC}}$ should only be used when the hypothesis of log-concavity can be expected to hold, at least approximately.
\begin{thm}
\label{Thm:Classifiers}
\begin{description}
\item[(a)] Assume $P_k \in \mathcal{P}_1$ for $k = 1,\ldots,K$.  Let $\mathcal{X}^* = \{x \in \mathbb{R}^d: |\argmax_k \pi_k f_k^*(x)| = 1\}$.  Then $\hat{C}_n^{\mathrm{LC}}(x) \stackrel{a.s.}{\rightarrow} C^{\mathrm{LCBayes}}(x)$ for almost all $x \in \mathcal{X}^*$, and
\[
\mathrm{Risk}(\hat{C}_n^{\mathrm{LC}}) \rightarrow \mathrm{Risk}(C^{\mathrm{LCBayes}}).
\]
\item[(b)] Now assume $P_k \in \mathcal{P}_2$ for $k = 1,\ldots,K$.  Let $\mathcal{X}^{**} = \{x \in \mathbb{R}^d: |\argmax_k \pi_k f_k^{**}(x)| = 1\}$.  Then $\hat{C}_n^{\mathrm{SLC}}(x) \stackrel{a.s.}{\rightarrow} C^{\mathrm{SLCBayes}}(x)$ for almost all $x \in \mathcal{X}^{**}$, and 
\[
\mathrm{Risk}(\hat{C}_n^{\mathrm{SLC}}) \rightarrow \mathrm{Risk}(C^{\mathrm{SLCBayes}}).
\]
\end{description}
\end{thm}
In fact, the smoothed log-concave classifier is somewhat easier to apply in practical classification problems than its unsmoothed analogue.  This is because if $x_0 \in \mathbb{R}^d$ is outside the convex hull of the training data for each of the $K$ classes (an event of positive probability), then the log-concave maximum likelihood estimates of the densities at $x_0$ are all zero.  Thus all such points would be assigned by $\hat{C}_n^{\mathrm{LC}}$ to Class 1.  On the other hand, $\hat{C}_n^{\mathrm{SLC}}$ avoids this problem altogether.  For these reasons, we considered only $\hat{C}_n^{\mathrm{SLC}}$ in our simulation study \citep{ChenSamworth2011} and below.

We remark that the direct use of $\hat{C}_n^{\mathrm{SLC}}$ (or any other classifier based on nonparametric density estimation) is not recommended when $d > 4$, due to the curse of dimensionality.  In such circumstances there are two options: dimension reduction (cf.~Section~\ref{Sec:BreastCancer} below), or further modelling assumptions such as independent component analysis models \citep{SamworthYuan2012}.  In either case, the methodology we develop remains applicable, but now as part of a more involved procedure.

\subsection{Breast cancer example}
\label{Sec:BreastCancer}
\begin{figure}[ht!]
 \centering
 $\begin{array}{c c}
  \includegraphics[scale=0.76]{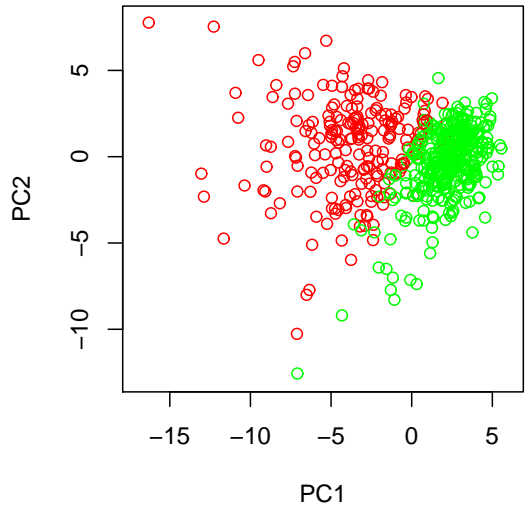} &
  \includegraphics[scale=0.30]{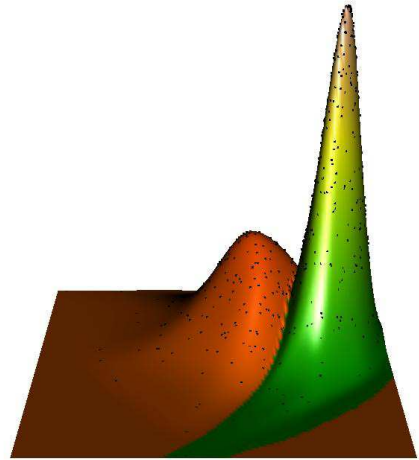} \\
  \mathrm{(a)} & \mathrm{(b)}             \\
  \includegraphics[scale=0.76]{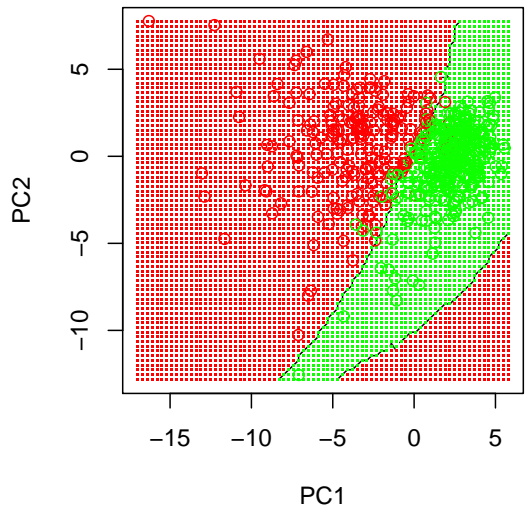} &
  \includegraphics[scale=0.76]{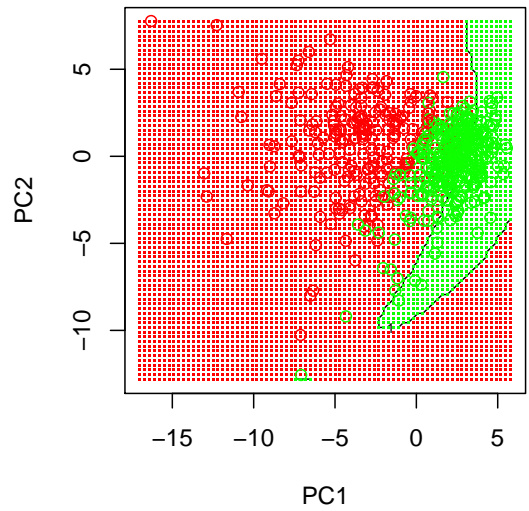} \\
  \mathrm{(c)} & \mathrm{(d)}
 \end{array}$
 \caption{(a) Wisconsin breast cancer data (benign cases in green; malignant
cases in red); (b) smoothed log-concave maximum likelihood density estimates; (c)
and (d) plot the decision boundaries of the smoothed log-concave classifier,
where the loss $L_2 = 1$ and $L_2=100$, respectively.}
 \label{Fig:wbcd}
\end{figure}

In the Wisconsin breast cancer data set \citep{SWM1993}, 30 measurements were taken from a digitised image of a fine needle aspirate of different breast masses.  There are 357 benign and 212 malignant instances, and we aim to construct a classifier based on this training data set to aid future diagnoses.  Only the first two principal components of the training data were considered, and these capture 63\% of the total variability; cf. Figure~\ref{Fig:wbcd}(a).  This was done to make our procedure computationally feasible, to reduce the effect of the curse of dimensionality, and to facilitate plots such as Figure~\ref{Fig:wbcd} below.

In Figure~\ref{Fig:wbcd}(b), we show the smoothed log-concave density estimates of both the benign and malignant classes. Figure~\ref{Fig:wbcd}(c) plots the decision boundaries of the smoothed log-concave classifier, where we treat benign cases and malignant cases equally.  However in practice, misdiagnosing a malignant tumour as benign is much more serious than misidentifying a benign one as malignant. One may therefore seek to incorporate different losses into the classifier. For $k = 1,2$, let $L_k$ denote the cost of failure to recognise the class $k$ (this notion can easily be generalised to multicategory situations were $L_{kk'}$ is the loss incurred in assigning the pair $(X,Y)$ to class $k'$ when $Y=k$).  Redefining the risk as 
\[
\mathrm{Risk}(C) = L_1\mathbb{P}\{C(X) = 2 \cap Y=1\} + L_2\mathbb{P}\{C(X) = 1 \cap Y=2\},
\]
the same asymptotic properties continue to hold, mutatis mutandis, for the classifier
\[
\hat{C}_n^{\mathrm{SLC*}}(x) = \argmax_{k \in \{1,2\}} N_k L_k \tilde{f}_{n,k}(x).
\]
We observe that this modification requires no recalculation of the smoothed log-concave density estimates and there is no loss of generality in taking $L_1 = 1$.  A GUI with slider is implemented in the \texttt{R} package \texttt{LogConcDEAD}, which provides a way of demonstrating how the decision boundaries change as $L_2$ varies. For the purpose of illustration, Figure~\ref{Fig:wbcd}(d) plots the decision boundaries of $\hat{C}_n^{\mathrm{SLC*}}$ when the cost $L_2$ of misidentifying a malignant tumour is 100. Compared with Figure~\ref{Fig:wbcd}(c), observations are of course considerably more likely to be classified as malignant under this setting. 

\subsection{Functional estimation problems}
\label{Sec:Functionals}

Classification problems are an important example of a situation where one is interested in a functional of one or more density estimates, rather than the density estimate itself.  For simplicity of exposition, we return in this section to the situation where we have a single independent sample $X_1,\ldots,X_n$ distributed according to a distribution $P_0$.  

In general, we can consider estimating a functional $\theta_0 = \theta(P_0)$ using the plug-in smoothed log-concave estimate $\tilde{\theta}_n = \theta(\tilde{P}_n)$, where $\tilde{P}_n$ is the distribution with density $\tilde{f}_n$.  Note that even if this functional cannot be computed directly, it is usually straightforward to construct a Monte Carlo approximation to $\tilde{\theta}_n$ by applying the algorithm for sampling from $\tilde{f}_n$ outlined in Section~\ref{Sec:Sampling}.  To describe the theoretical properties of these functional estimates, for $a > 0$, let $\mathcal{B}_a$ denote the set of signed measures $P$ on $\mathbb{R}^d$ with $\int_{\mathbb{R}^d} e^{a\|x\|} \, d|P|(x) < \infty$.  Equip $\mathcal{B}_a$ with the norm
\[
\|P\|_a = \int_{\mathbb{R}^d} e^{a\|x\|} \, d|P|(x). 
\]
We can then consider $\theta$ as a measurable function on $(\mathcal{B}_a,\|\cdot\|_a)$ taking values in some other normed space $(\mathcal{B},\|\cdot\|)$.  
\begin{prop}
\label{Prop:Functionals}
Let $P_0 \in \mathcal{P}_2$, and let $P_0^{**}$ denote the probability distribution whose density is the smoothed version of the log-concave approximation to $P_0$.  Suppose that $\theta:\mathcal{B}_a \rightarrow \mathcal{B}$ is continuous, and let $\theta^{**} = \theta(P_0^{**})$.  Then $\|\tilde{\theta}_n - \theta^{**}\| \stackrel{a.s.}{\rightarrow} 0$ as $n \rightarrow \infty$.
\end{prop}
Once again, we remark that if $P_0$ has a log-concave density, then $P_0 = P_0^{**}$.  The fact that the topology on $\mathcal{B}_a$ is rather strong means that the continuity requirement on $\theta$ is relatively weak.  This is illustrated in the following corollary, which considers the special case of \emph{linear} functionals in Proposition~\ref{Prop:Functionals}.
\begin{cor}
\label{Cor:LinearFunctionals}
Let $P_0 \in \mathcal{P}_2$, and let $a_0 > 0$ and $b_0 \in \mathbb{R}$ be such that $f^{**}(x) \leq e^{-a_0\|x\| + b_0}$, where $f^{**}$ is the smoothed log-concave approximation to $P_0$.  Let $\theta(P) = \int_{\mathbb{R}^d} g \, dP$ for some measurable function $g:\mathbb{R}^d \rightarrow \mathbb{R}$ satisfying
\begin{equation}
\label{Eq:FuncCond}
\sup_{x \in \mathbb{R}^d} e^{-a\|x\|}|g(x)| < \infty
\end{equation}
for some $a < a_0$.  Then $\tilde{\theta}_n \stackrel{a.s.}{\rightarrow} \theta^{**}$.
\end{cor}

\section*{Acknowledgments} We thank the Associate Editor and two anonymous referees for their helpful comments.  The second author is grateful for the support of a Leverhulme Research Fellowship and an EPSRC Early Career Fellowship.

\section{Appendix}

\begin{prooftitle}{of Proposition~\ref{Prop:Basic}}

\textbf{(a)} This follows immediately from Theorems~2.8 and~2.18 of \citet{DharmadhikariJoag-Dev1988}.

\textbf{(b)} Note that for any non-empty open set $B \subseteq \mathbb{R}^d$,
\[
\int_B \tilde{f}_n(x) \, dx = \int_B \int_{C_n} \hat{f}_n(y) \phi_{d,\hat{A}}(x-y) \, dy \, dx,
\]
which is positive, since the integrand is positive and continuous on the region of integration.

\textbf{(c)} The fact that $\tilde{f}_n$ is infinitely differentiable follows from Proposition~8.10 of \citet{Folland1999}.  In fact, using standard multi-index notation with $\alpha = (\alpha_1,\ldots,\alpha_d)$ and $\partial^\alpha = \bigl(\frac{\partial}{\partial x_1}\bigr)^{\alpha_1}\ldots \bigl(\frac{\partial}{\partial x_1}\bigr)^{\alpha_d}$, we have $\partial^\alpha \tilde{f}_n = \hat{f}_n \ast \partial^\alpha \phi_{d,\hat{A}}$.  Writing $|\alpha| = \sum_{l=1}^d \alpha_l$ and $\alpha ! = \prod_{l=1}^d \alpha_l!$, it follows that for any $x_0 \in \mathbb{R}^d$ and $k \in \mathbb{N}$, 
\begin{align*}
\biggl|\tilde{f}_n(x) - \sum_{|\alpha| \leq k} \frac{(\partial^\alpha \tilde{f}_n)(x_0)}{\alpha !}(x\!-\!x_0)^\alpha\biggr| &\leq \int_{C_n} \hat{f}_n(y)\biggl|\phi_{d,\hat{A}}(x\!-\!y) - \sum_{|\alpha| \leq k} \frac{(\partial^\alpha \phi_{d,\hat{A}})(x_0)}{\alpha !}(x\!-\!y\!-\!x_0)^\alpha\biggr| \, dy \\
&\rightarrow 0
\end{align*}
as $k \rightarrow \infty$, by the dominated convergence theorem and Lemma~1 of \citet{CuleSamworth2010}.

\textbf{(d)} Conditional on $X_1,\ldots,X_n$, let $X^*$ and $Y^*$ be independent, with $X^*$ having density $\hat{f}_n$ and $Y^*$ having density $\phi_{d,\hat{A}}$, so that $X^* + Y^*$ has conditional density $\tilde{f}_n$.  Then
\[
\mathbb{E}(X^*+Y^*|X_1,\ldots,X_n) =  \mathbb{E}(X^*|X_1,\ldots,X_n) = \int_{\mathbb{R}^d} x \hat{f}_n(x) \, dx = \bar{X},
\]
and
\[
\mathrm{Cov}(X^*+Y^*|X_1,\ldots,X_n) = \mathrm{Cov}(X^*|X_1,\ldots,X_n) + \mathrm{Cov}(Y^*|X_1,\ldots,X_n) = \tilde{\Sigma} + {\hat{A}} = \hat{\Sigma}.
\]
\hfill $\Box$
\end{prooftitle}

\begin{prooftitle}{of Theorem~\ref{Thm:Asymp}}
Let $d_P$ and $d_{TV}$ denote the Prohorov and total variation metrics on the space of probability measures on $\mathbb{R}^d$.  Recall that $d_P$ metrises weak convergence, and that $d_P \leq d_{TV}$.  Let $\hat{\mu}_n$ denote the probability measure corresponding to the density $\tilde{f}_n$, let $\hat{\nu}_n$ denote the probability measure corresponding to the convolution of $f^*$ with the measure $N_d(0,\hat{A})$, and let $\nu$ denote the probability measure corresponding to $f^{**}$.  Then
\begin{align}
\label{Eq:TwoTerms}
d_P(\hat{\mu}_n,\nu) &\leq d_P(\hat{\mu}_n,\hat{\nu}_n) + d_P(\hat{\nu}_n,\nu) \nonumber \\
&\leq d_{TV}(\hat{\mu}_n,\hat{\nu}_n) + d_P(\hat{\nu}_n,\nu) \nonumber \\
&= \frac{1}{2}\int_{\mathbb{R}^d} |\hat{f}_n \ast N_d(0,\hat{A}) - f^* \ast N_d(0,\hat{A})| + d_P(\hat{\nu}_n,\nu) \nonumber \\
&\leq \frac{1}{2} \int_{\mathbb{R}^d} |\hat{f}_n - f^*| + d_P(\hat{\nu}_n,\nu).
\end{align}
The first term of~(\ref{Eq:TwoTerms}) converges almost surely to zero, by Theorem~2.15 of \citet{DSS2011}.  The second term also converges almost surely to zero, using the fact that $\hat{A} \stackrel{a.s.}{\rightarrow} A^*$ as $n \rightarrow \infty$.  Proposition~2 of \citet{CuleSamworth2010} strengthens the mode of convergence and yields the result.
\hfill $\Box$
\end{prooftitle}

\begin{prooftitle}{of Proposition~\ref{Prop:Bounds}}
If $x \in C_{n,j}$, and $\hat{f}_n(x) = \exp(b_j^T x - \beta_j)$, then $\hat{f}_n(x-y) \leq \exp\{b_j^T (x-y) - \beta_j\}$ for all $y \in \mathbb{R}^d$.  It follows that 
\begin{equation}
\label{Eq:InternalUpper}
\frac{\tilde{f}_n(x) - \hat{f}_n(x)}{\hat{f}_n(x)} \leq \int_{\mathbb{R}^d} e^{-b_j^T y}\phi_{d,\hat{A}}(y) \, dy - 1 = e^{\frac{1}{2}b_j^T \hat{A}b_j} - 1.
\end{equation}
Now
\begin{align*}
\int_{\mathbb{R}^d} |\tilde{f}_n - \hat{f}_n| &= \int_{C_n} |\tilde{f}_n - \hat{f}_n| + \delta_n \\
&= \int_{C_n} (\tilde{f}_n - \hat{f}_n)_+ + \int_{C_n} (\hat{f}_n - \tilde{f}_n)_+ + \delta_n.
\end{align*}
But
\[
\int_{C_n} (\hat{f}_n - \tilde{f}_n)_+ = \int_{C_n} (\tilde{f}_n - \hat{f}_n)_+ - \int_{C_n} (\tilde{f}_n - \hat{f}_n) = \int_{C_n} (\tilde{f}_n - \hat{f}_n)_+ + \delta_n.
\]
It therefore follows from this and~(\ref{Eq:InternalUpper}) that 
\[
\int_{\mathbb{R}^d} |\tilde{f}_n - \hat{f}_n| \leq 2\sum_{j \in J} \int_{C_{n,j}} \hat{f}_n(x) (e^{\frac{1}{2}b_j^T \hat{A}b_j} - 1) \, dx + 2 \delta_n \leq 2(e^{\frac{1}{2}\lambda_{\max}} - 1 + \delta_n),
\]
as required.
\hfill $\Box$
\end{prooftitle}

\begin{prooftitle}{of Theorem~\ref{Thm:Independent}}
\textbf{(a)} Let $f$ be an arbitrary log-concave density on $\mathbb{R}^d$, and let $X$ be a random vector with density $f$.  Letting $X = (X_1^T,X_2^T)^T$, where $X_1$ and $X_2$ take values in $\mathbb{R}^{d_1}$ and $\mathbb{R}^{d_2}$ respectively, we write $f_{X_1}$ for the marginal density of $X_1$ and $f_{X_2|X_1}(\cdot|x_1)$ for the conditional density of $X_2$ given $X_1=x_1$.  By Theorem~6 of \citet{Prekopa1973}, $f_{X_1}$ is log-concave and by Proposition~1 of \citet{CSS2010}, $f_{X_2|X_1}(\cdot|x_1)$ is log-concave for each $x_1$.  

There is also no loss of generality in assuming $f$ is upper semi-continuous.  Since $P \in \mathcal{P}_1$, we may assume without loss of generality that $\int_{\mathbb{R}^d} |\log f| \, dP < \infty$.  We may therefore apply Fubini's theorem and seek to maximise over all upper semi-continuous log-concave densities the quantity
\begin{align}
\label{Eq:Fubini}
\int_{\mathbb{R}^d} \log f \, dP &= \int_{\mathbb{R}^{d_1}} \int_{\mathbb{R}^{d_2}} \{\log f_{X_1}(x_1) + \log f_{X_2|X_1}(x_2|x_1)\} dP_2(x_2) dP_1(x_1) \nonumber \\
&= \int_{\mathbb{R}^{d_1}} \log f_{X_1}(x_1) P_1(dx_1) + \int_{\mathbb{R}^{d_1}} \int_{\mathbb{R}^{d_2}} \log f_{X_2|X_1}(x_2|x_1) dP_2(x_2) dP_1(x_1).
\end{align}
The first term on the right-hand side of~(\ref{Eq:Fubini}) is maximised uniquely over all upper semi-continuous log-concave densities by setting $f_{X_1} = f_1^*$.  Moreover, for any fixed $x_1$, the quantity $\int_{\mathbb{R}^{d_2}} \log f_{X_2|X_1}(x_2|x_1) \, dP_2(x_2)$ is maximised uniquely over upper semi-continuous log-concave densities by setting $f_{X_2|X_1}(\cdot|x_1) = f_2^*$.  Since this choice does not depend on $x_1$, it maximises the second term on the right-hand side of~(\ref{Eq:Fubini}).  Because both terms can be maximised simultaneously, it follows that $f^* = f_1^*f_2^*$, as desired.

\textbf{(b)} Write $\Sigma$ and $\Sigma^*$ for the covariance matrices corresponding to the probability distribution $P$ and the density $f^*$ respectively.  The independence structure of $P_0$ and $f^*$ gives that 
$\Sigma = \left[ \begin{array}{cc}
\Sigma_1 & 0 \\
0 & \Sigma_2 \end{array} \right]$ and
$\Sigma^* = \left[ \begin{array}{cc}
\Sigma^*_1 & 0 \\
0 & \Sigma^*_2 \end{array} \right]$.  Here, $\Sigma_1$ and $\Sigma_1^*$ are $d_1 \times d_1$ submatrices, while $\Sigma_2$ and $\Sigma_2^*$ are $d_2 \times d_2$ submatrices.  Therefore, $A^* = \Sigma - \Sigma^*$ is of the form 
$A^* = \left[ \begin{array}{cc}
A^*_1 & 0 \\
0 & A^*_2 \end{array} \right]$.  Writing $x,y \in \mathbb{R}^d$ as $(x_1^T,x_2^T)^T$ and $(y_1^T,y_2^T)^T$ respectively, where $x_1,y_1 \in \mathbb{R}^{d_1}$ and $x_2,y_2 \in \mathbb{R}^{d_2}$, it follows again by Fubini's theorem that
\begin{align*}
f^{**}(x) &= (f^* \ast N_d(0,A^*))(x) \\
&= \int_{\mathbb{R}^{d_1}} \int_{\mathbb{R}^{d_2}} f_1^*(y_1) f_2^*(y_2) \, dN_{d_2}(0,A_2^*)(x_2-y_2) \, dN_{d_1}(0,A_1^*)(x_1-y_1) \\
&=\biggl\{\int_{\mathbb{R}^{d_1}} f_1^*(y_1) \, dN_{d_1}(0,A_1^*)(x_1-y_1)\biggr\}\biggl\{\int_{\mathbb{R}^{d_2}} f_2^*(y_2) \, dN_{d_2}(0,A_2^*)(x_2-y_2) \biggr\} \\
&= f^{**}_1(x_1) f^{**}_2(x_2).
\end{align*}
\hfill $\Box$
\end{prooftitle}


\begin{prooftitle}{of Theorem~\ref{Thm:Covariance}}
Let $P \in \mathcal{P}_1$, and let $f^*$ denote its log-concave approximation.  Without loss of generality, we may assume $\int_{\mathbb{R}^d} x \, dP(x) = 0$, so it suffices to show that if $A^* := \int_{\mathbb{R}^d} xx^T \, dP(x) - \int_{\mathbb{R}^d} xx^T f^*(x) \, dx$ is the zero matrix, then $P$ has a log-concave density.

Let $P^*$ denote the distribution corresponding to $f^*$, let $X \sim P$ and let $X^* \sim P^*$.  For an arbitrary $u \in \mathbb{R}^d$, let $F_u$ and $F_u^*$ denote the distribution functions of $u^T X$ and $u^T X^*$ respectively, and let
\[
	G_u(s)  =  \int_{-\infty}^s \, F_u(t) \, dt \quad \mbox{ and } \quad 
	G^*_u(s)  =  \int_{-\infty}^s \, F_u^*(t) \, dt.
\]
Fix $s \in \mathbb{R}$.  By applying Remark~2.3 of \citet{DSS2011} to the convex function $x \mapsto (s - u^T x)_+$ and Fubini's theorem, we have that 
\begin{align}
\nonumber
	0 \leq \int_{\mathbb{R}^d} (s- u^T x)_+ \, d(P-P^*)(x) &= \int_{\mathbb{R}^d} \int_{-\infty}^\infty \mathbbm{1}_{\{ u^T x \leq t < s\}} \, dt \, d(P-P^*)(x) \\
\label{Eq:cdfmarginal}
          &= \int_{-\infty}^{s} (F_u-F_u^*)(t) \, dt  =  G_u(s)-G_u^*(s).
\end{align}
Since all moments of log-concave densities are finite, we have $\int_{\mathbb{R}^d} xx^T f^*(x) \, dx < \infty$.  So, since $A^* = 0$, we must have $P \in \mathcal{P}_2$.  We can therefore integrate by parts as follows:
\begin{align}
0 = \int_{\mathbb{R}^d} (u^T x)^2 \, d(P-P^*)(x) = \int_{-\infty}^\infty t^2 \, d(F_u-F_u^*)(t) &= - 2 \int_{-\infty}^{\infty} \, t (F_u-F_u^*)(t) \, dt \nonumber \\ \label{Eq:cdfmarginal2}
&= 2 \int_{-\infty}^{\infty} (G_u-G_u^*)(t) \, dt. 
\end{align}
Combining~(\ref{Eq:cdfmarginal}),~(\ref{Eq:cdfmarginal2}) and the fact that $G_u - G_u^*$ is continuous, we deduce that $G_u=G_u^*$.   Thus $F_u=F_u^*$, by the fundamental theorem of calculus and the fact that $F_u$ and $F_u^*$ are both right-continuous.  It follows that 
\[
\mathbb{E}(e^{iu^TX}) = \int_{-\infty}^\infty e^{it} \, dF_u(t) = \int_{-\infty}^\infty e^{it} \, dF_u^*(t) = \mathbb{E}(e^{iu^TX^*}).
\]
Since $u \in \mathbb{R}^d$ was arbitrary, we deduce that $P=P^*$, so $P$ has a log-concave density. 
\hfill $\Box$
\end{prooftitle}

\begin{prooftitle}{of Proposition~\ref{Prop:Convex}}
Suppose that the upper semi-continuous log-concave density $f^*$ is the log-concave approximation to $P_1,P_2 \in \mathcal{P}_1$.  Then for each $t \in (0,1)$, we see that $f^*$ also maximises
\[
\int_{\mathbb{R}^d} \log f \, d(tP_1 + (1-t)P_2) = t \int_{\mathbb{R}^d} \log f \, dP_1 + (1-t) \int_{\mathbb{R}^d} \log f \, dP_2
\]
over all upper semi-continuous log-concave densities $f$ on $\mathbb{R}^d$.
\hfill $\Box$
\end{prooftitle}

\begin{prooftitle}{of Theorem~\ref{Thm:Test}}
Let $d_2$ denote the second Mallows metric on $\mathcal{P}_2$, so $d_2(P,Q) = \inf_{(X,Y) \sim (P,Q)} \{\mathbb{E}\|X-Y\|^2\}^{1/2}$, where the infimum is taken over all pairs $(X,Y)$ of random vectors $X \sim P$ and $Y \sim Q$ on a common probability space.  Recall that the infimum in this definition is attained, and that if $P,P_1,P_2,\ldots \in \mathcal{P}_2$, then $d_2(P_n,P) \rightarrow 0$ if and only if both $P_n \stackrel{d}{\rightarrow} P$ and $\int_{\mathbb{R}^d} \|x\|^2 \, dP_n(x) \rightarrow \int_{\mathbb{R}^d} \|x\|^2 \, dP(x)$.  Let $P^*$ denote the distribution corresponding to the log-concave approximation to $P_0$, and for $\delta > 0$ to be chosen later, let $\mathcal{Q}_{2,\delta}$ denote the subset of $\mathcal{P}_2$ consisting of those distributions $Q$ with $d_2(Q,P^*) \leq \delta$ that have a log-concave density.  Fix $\epsilon > 0$ and let $Q \in \mathcal{Q}_{2,\delta}$.  Let $\mathbb{P}_n$ and $\mathbb{Q}_n$ denote the empirical distribution of an independent sample of size $n$ from $P^*$ and an independent sample from $Q$ respectively.  We will require a bound for $\mathbb{P}\{d_2(\mathbb{Q}_n,\mathbb{P}_n) \geq \epsilon/4\}$ that holds uniformly over $\mathcal{Q}_{2,\delta}$, and obtain this using the following coupling argument.  We may suppose that $(X_1,Y_1),\ldots,(X_n,Y_n)$ are independent and identically distributed pairs with $X_i \sim P^*$ and $Y_i \sim Q$ and that $\mathbb{P}_n$ and $\mathbb{Q}_n$ are obtained as the empirical distribution of $X_1,\ldots,X_n$ and $Y_1,\ldots,Y_n$ respectively.  We may further suppose that $\mathbb{E}\|X_i - Y_i\|^2 = d_2^2(P^*,Q)$; in other words, $X_i$ and $Y_i$ are coupled in such a way that they attain the infimum in the definition of the second Mallows distance.  
Using standard results on the Mallows distance (e.g. Equation~(8.2) and Lemma~8.7 of \citet{BickelFreedman1981}), we deduce that for $\delta \leq \epsilon^{3/2}/32$,
\begin{align*}
\sup_{Q \in \mathcal{Q}_{2,\delta}} \mathbb{P}\{d_2(\mathbb{Q}_n,\mathbb{P}_n) > \epsilon/4\} &\leq \sup_{Q \in \mathcal{Q}_{2,\delta}} \mathbb{P}\biggl(\frac{1}{n}\sum_{i=1}^n \|X_i - Y_i\|^2 > \frac{\epsilon^2}{16}\biggr) \\
&\leq \frac{16}{\epsilon^2} \sup_{Q \in \mathcal{Q}_{2,\delta}} \mathbb{E}(\|X_1 - Y_1\|^2) \leq  \frac{16\delta^2}{\epsilon^2} \leq \frac{\epsilon}{2}.
\end{align*}
Now let $\hat{Q}_n$ denote the distribution corresponding to the log-concave maximum likelihood estimator constructed from $X_1,\ldots,X_n$, and let $\hat{\mathbb{Q}}_n$ denote the empirical distribution of a sample of size $n$ which, conditional on $X_1,\ldots,X_n$, is drawn independently from $\hat{Q}_n$.  By reducing $\delta > 0$ if necessary, we may assume $\delta \leq \epsilon/4$.  It follows that
\begin{align}
\label{Eq:Mallows}
\mathbb{P}\{d_2(\hat{\mathbb{Q}}_n,P^*) > \epsilon\} &\leq \sup_{Q \in \mathcal{Q}_{2,\delta}} \mathbb{P}\{d_2(\mathbb{Q}_n,Q) > 3\epsilon/4\} + \mathbb{P}\{d_2(\hat{Q}_n,P^*) > \delta\} \nonumber \\
&\leq \sup_{Q \in \mathcal{Q}_{2,\delta}} \mathbb{P}\{d_2(\mathbb{Q}_n,\mathbb{P}_n) > \epsilon/4\} + \mathbb{P}\{d_2(\mathbb{P}_n,P^*) > \epsilon/4\} + \mathbb{P}\{d_2(\hat{Q}_n,P^*) > \delta\} \nonumber \\
&\leq \frac{\epsilon}{2} + \mathbb{P}\{d_2(\mathbb{P}_n,P^*) > \epsilon/4\} + \mathbb{P}\{d_2(\hat{Q}_n,P^*) > \delta\} \leq \epsilon
\end{align}
for sufficiently large $n$.  The final convergence of the second term here follows from the weak law of large numbers, while for the third term it follows from Proposition~2(c) of \citet{CuleSamworth2010} and the dominated convergence theorem.

Let $\hat{\mathbb{Q}}_{nb}$ and $\hat{Q}_{nb}$ denote respectively the empirical distribution and the distribution corresponding to the log-concave maximum likelihood estimator of the $b$th bootstrap sample $X_{1b}^*,\ldots,X_{nb}^*$ drawn from $\hat{Q}_n$.  We deduce from~(\ref{Eq:Mallows}), Theorem 2.15 of \citet{DSS2011} and another application of Proposition~2(c) of \citet{CuleSamworth2010} that there exists $a > 0$ such that 
\begin{equation}
\label{Eq:hatQ}
\int_{\mathbb{R}^d} e^{a\|x\|} \, d(\hat{Q}_{nb} - P^*)(x) \stackrel{p}{\rightarrow} 0.
\end{equation}
Now let 
\[
\hat{A}_{nb} \equiv \hat{\Sigma}_b - \tilde{\Sigma}_b \equiv  \frac{n}{n-1} \int_{\mathbb{R}^d} (x - \bar{X}_b^*)(x - \bar{X}_b^*)^T \, d\hat{\mathbb{Q}}_{nb}(x) - \int_{\mathbb{R}^d} (x - \bar{X}_b^*)(x - \bar{X}_b^*)^T \, d\hat{Q}_{nb}(x),
\]
where $\bar{X}_b^* = n^{-1}\sum_{i=1}^n X_{ib}^*$.  From~(\ref{Eq:Mallows}),~(\ref{Eq:hatQ}), the dominated convergence theorem and the continuous mapping theorem, we have that $\operatorname{tr}(\hat{A}_{nb}) \stackrel{p}{\rightarrow} 0$ as $n \rightarrow \infty$.  On the other hand, in the notation of Theorem~\ref{Thm:Asymp},
\[
\operatorname{tr}(\hat{A}) = \operatorname{tr}(\hat{\Sigma}) - \operatorname{tr}(\tilde{\Sigma}) \stackrel{p}{\rightarrow} \operatorname{tr}(\Sigma) - \operatorname{tr}(\Sigma^*) = \operatorname{tr}(A^*)> 0,
\]
where the final claim follows from Theorem~\ref{Thm:Covariance} and the fact that $P_0$ does not have a log-concave density.  Note that this claim holds even if $P_0 \in \mathcal{P}_1 \setminus \mathcal{P}_2$, in which case $\operatorname{tr}(\Sigma) = \infty$.

Write $Z_{nb} = \mathbbm{1}_{\{\operatorname{tr}(\hat{A}_{nb}) > \operatorname{tr}(A^*)/2\}}$, and note that $Z_{n1},\ldots,Z_{nB}$ are exchangeable (so in particular, identically distributed).  Thus, for any $\alpha \in (0,1)$,
\begin{align*}
\mathbb{P}(\text{Do not reject} \ H_0) &= \mathbb{P}\biggl(\frac{1}{B+1}\sum_{b=1}^{B+1} \mathbbm{1}_{\{\operatorname{tr}(\hat{A}) > \operatorname{tr}(\hat{A}_{nb})\}} \leq 1-\alpha\biggr) \\
&\leq \mathbb{P}\{\operatorname{tr}(\hat{A}) \leq \operatorname{tr}(A^*)/2\} + \mathbb{P}\biggl(\frac{1}{B+1}\sum_{b=1}^{B+1} Z_{nb} \geq 1 - \alpha\biggr) \\
&\leq \mathbb{P}\{\operatorname{tr}(\hat{A}) \leq \operatorname{tr}(A^*)/2\} + \frac{\mathbb{E}(Z_{n1})}{1-\alpha} \rightarrow 0
\end{align*}
as $n \rightarrow \infty$.  We deduce that for any given size of test $\alpha \in (0,1)$, the power at any alternative converges to 1.
\hfill $\Box$
\end{prooftitle}

\begin{prooftitle}{of Theorem~\ref{Thm:Classifiers}}
\textbf{(a)} Note that  
\[
\hat{C}_n^{\mathrm{LC}}(x) = \argmax_{k \in \{1,\ldots,K\}} \frac{N_k}{n} \hat{f}_{n,k}(x).
\]
We have that $\int_{\mathbb{R}^d} |\hat{f}_{n,k} - f^*_k| \stackrel{a.s.}{\rightarrow} 0$ as $n \rightarrow \infty$ for every $k$, and in fact, by Theorem~10.8 of \citet{Rockafellar1997}, it is almost surely the case that $\hat{f}_{n,k}$ converges to $f_k^*$ uniformly on compact sets in the interior of the support of $f_k^*$.  By the strong law of large numbers and the fact that the boundary of the support of $f_k^*$ has zero $d$-dimensional Lebesgue measure, it therefore follows that  
\[
\hat{C}_n^{\mathrm{LC}}(x) \stackrel{a.s.}{\rightarrow} C^{\mathrm{LCBayes}}(x)
\]
for almost all $x \in \mathcal{X}^*$.

In fact, with probability one, $\frac{N_k}{n}\hat{f}_{n,k}$ converges to $\pi_k f_k^*$ uniformly on compact sets in the interior of the support of $f_k^*$.  It follows immediately from this and the dominated convergence theorem that 
\[
\mathrm{Risk}(\hat{C}_n^{\mathrm{LC}}) \rightarrow \mathrm{Risk}(C^{\mathrm{LCBayes}}).
\]

\textbf{(b)} The proof is virtually identical to that of Part (a), so is omitted.
\hfill $\Box$
\end{prooftitle}

\begin{prooftitle}{of Proposition~\ref{Prop:Functionals}}
The conclusion of Theorem~\ref{Thm:Asymp} can be stated in the notation of Section~\ref{Sec:Functionals} as
\[
\|\tilde{P}_n - P_0^{**}\|_a \stackrel{a.s.}{\rightarrow} 0.
\]
The result therefore follows immediately by the continuous mapping theorem.
\hfill $\Box$
\end{prooftitle}

\begin{prooftitle}{of Corollary~\ref{Cor:LinearFunctionals}}
It suffices to show that under condition~(\ref{Eq:FuncCond}), the functional $\theta(P) = \int_{\mathbb{R}^d} g \, dP$ is continuous.  Fix $a < a_0$ such that $\sup_{x \in \mathbb{R}^d} e^{-a\|x\|}|g(x)| < \infty$, and choose a sequence $(P_n)$ such that $\|P_n - P\|_a \rightarrow 0$.  Then
\begin{align*}
|\theta(P_n) - \theta(P)| &\leq \int_{\mathbb{R}^d} |g| \, d|P_n-P| \\
&\leq \sup_{x \in \mathbb{R}^d} e^{-a\|x\|}|g(x)|\int_{\mathbb{R}^d} e^{a\|x\|} \, d|P_n - P| \\
&= \sup_{x \in \mathbb{R}^d} e^{-a\|x\|}|g(x)| \|P_n - P\|_a \rightarrow 0
\end{align*}
as $n \rightarrow \infty$.  Thus $\theta$ is continuous, as required.
\hfill $\Box$
\end{prooftitle}

\end{document}